\documentclass{article}

\usepackage{latexsym}
\usepackage{amsmath}
\usepackage{amsfonts}
\usepackage{graphicx}
\usepackage{color}
\usepackage{hyperref}
\newcommand{\ds}{\displaystyle}

\newcommand{\torus}{{\mathbb{T}}}
\newcommand{\gu}{\mathbf{u}}

\usepackage{theorem}

\theoremheaderfont{\scshape}
\newtheorem{theorem}{Theorem}[section]

\newtheorem{prop1}[theorem]{Proposition}
\theorembodyfont{\normalfont}

\newenvironment{proof}{{\bf Proof }}{\hbox{~} \hfill \rule{0.5em}{0.5em}\\}
\numberwithin{equation}{section}

\begin{document}

\title{Homogenization and transport equations: the case of desert and sand piles\footnote{{This work is supported by  NLAGA Project (Non Linear Analysis, Geometry and Applications Project).}}}         
\date{}          

\maketitle

\centerline{\scshape Badahi ould Mohamed \footnote{badahi1977@yahoo.fr }}
\medskip
{\footnotesize
 \centerline{ Al-Jouf University, Kingdom of Saudi Arabia.}
\centerline{  Faculty of Science and Arts  at Al Qurayat }
   \centerline{Department of Mathematics}
} 
\medskip
\centerline{\scshape Ibrahima Faye \footnote{ibrahima.faye@uadb.edu.sn }}
\medskip
{\footnotesize
 \centerline{Universit\'e Alioune Diop de Bambey, UFR S.A.T.I.C, BP 30 Bambey (S\'en\'egal),}
   \centerline{Laboratoire de Math\'ematiques de la D\'ecision et d'Analyse Num\'erique}
\centerline{ (L.M.D.A.N). }
} 

\medskip
\centerline{\scshape Diaraf Seck \footnote{diaraf.seck@ucad.edu.sn}}
\medskip
{\footnotesize
 \centerline{Universit\'e Cheikh Anta Diop de Dakar, BP 16889 Dakar Fann,}
  \centerline{  Ecole Doctorale de
                      Math\'ematiques et Informatique. }
   \centerline{Laboratoire de Math\'ematiques de la D\'ecision et d'Analyse Num\'erique}
    \centerline{ (L.M.D.A.N) F.A.S.E.G)/F.S.T. }

}
\pagestyle{myheadings}
 \renewcommand{\sectionmark}[1]{\markboth{#1}{}}
\renewcommand{\sectionmark}[1]{\markright{\thesection\ #1}}
\begin{abstract}\noindent In this paper we build models for short-term, mean-term and long-term
dynamics of dune in desert. They are models that
are degenerated parabolic equations which are, moreover, singularly perturbed.
We, then give existence and uniqueness results for the  models, followed by homogenization ones and a corrector result is given.
\end{abstract}
keywords: transport equations,dunes, PDE, homogenizaton, two scale convergence, corrector result.
\section{Introduction}
This paper deals about sand transport problems in desert. Desertifiction is a natural phenomenon caused by climate variations and now days by human activities.
The problem of aeolian transport of sand, which aroused from
few years many research represents interest from the
physicians community  of granular media. The aim ultimately is to achieve a comprehensive understanding of sand transport by wind, in order to understand the phenomena of formation and migration of the dunes in desert environments. This represents a challenge scientifically. First, the wind flow in the atmosphere is unsteady, turbulent and influenced by the earth's topography. The modeling of turbulent flows still raises many questions.
The objective will be to study mathematically, deformation and spread sand dunes under the flow effect of the air.
First we will explain how to build models with participation of parameters and framework. Then we plan to do a study
theoretical mathematical models proposed. This mathematical study will focus on the homogenization of tools by considering a periodic environment.

~
The location of the dune forms is directly dependent on the particle size of the soil particles. The wind exerts its action on well-defined size of materials.

Dune formation is related to the movement of sand particles. There are three different ways of driving the particles: saltation, creeping on the surface and the suspension.

- Saltation:

The initial movement of soil particles is a series of jumps. The particle diameter saltation is between 0.5 and 1.1 mm. After jumping, the particles fall under the action of gravity. The descending portion of the path is inclined towards the ground and substantially straight. Few particles reach an altitude of over 1 m and about 90 per cent of them are less than jumps 30 cm. The horizontal amplitude of a jump is generally between 0.5 and 1 m.\\
The saltation phenomenon is essential to initiate wind erosion. It is because of the two other modes of soil components by the wind creeping surface and airborne.

- The creep surface

Particles of larger size roll or slide overground. Too heavy to be lifted, their movement is triggered by the impact of particles saltation rather than by the action of the wind. The particles which move and have diameters between 0.5 and 2 mm depending on their density and wind speed.

- Suspension

Generally fine dust can be swept away as if they were projected into the air by the impact of larger grains. Once arrived in the turbulent layer they can be raised to great heights by air updrafts and form dust clouds often reaching heights of 3-4000 meters. Although their appearance may be impressive, the essential mechanism of wind erosion remains saltation because without such clouds could occur. For additional details see \cite{AndClauDou, Bagnold1, CokeWarreGou, Hersen}.\\
The paper is organized as follows: In the second section we give
the modeling of sand transport problems. In the section 3, we give existence and uniqueness result for the models. The section 4 deals with homogenized method and a corrector result is also given in the last section.

%
%
\section{Models of interest}
In this section, we present several models of linear parabolic equations described phenomena of transport. These models derived from the
Exner equation, which models sand transport. It is given by the expression
\begin{equation}\label{eq001}
\frac{\partial z}{\partial t}+\frac{1}{1-p}\nabla\cdot q=0,
\end{equation}
where $p$ is the porosity of the surface i.e. the percentage of void in the sediment,
$q$ is  the sand volume flow and $z = z (x, t)$ is the height of
the layer of sediment in position $x$ and at time $t.$\\

The flow of materials transported consider  the flow $q_r$ rated thrust and the quantity being transported in suspenion $q_s.$
There are two approaches to determine $q:$ the first one is to estimate separately the flows $q_r$ and $q_s$ and to sum to obtain the total flux and the second one is to estimate directly the total flow.

 \begin{itemize}
 \item The flow from the substances carried in suspension $ q_{s} $  expressed in $ m^{2}.s^{-1} $ by :\\
 $$ q_{s}=\int^{\eta}_{y_{0}} u_{p}(x,y,t) \phi(x,y,t)dy, $$
 with $ \eta$ being the free surface of the fluid ,$ y_{0} $ the thickness of the movable bottom layer,
$  u_{p}(x,y,t) $ is the speed of a grain and $ \phi(x,y,t) $  is the sediment concentration \\
\item The flow from the materials transported by thrusting: 
We will mention some approaches\\
b-1.\,\, Du Bays  pioneering approach : $ q_{r}=C_{B}F_{e}(F_{e}-F_{0}) $\\
where $ F_{e} $ is the shear stress ,$ C_{B} $ stands for the characteristic coefficient of  the flow and $F_{0}$ depends on layer thickness \\
b-2.\,\, The Yang approach: $ q_{r}=\frac{k}{C^{4}}\rho gU^{2}(U^{2}-U^{2}_{s}) $\\
where $ U $ is the flow velocity, $g$ is the gravity, $\rho$ is the density, $U_s$ is the flow of power and $C$ is the friction coefficient.
\end{itemize}
 So we shall mention a few of them according to the flow and friction velocity. These approaches are based on the Saint-Venant and Navier
Stokes equations.\\
\begin{enumerate}
\item

The first is related to the speed and is proposed by Gekerma \cite{Gek}. It is given by
\begin{equation}\label{eq2}
q=\alpha|\mathbf  u|^3\frac{\mathbf u}{|\mathbf u|}-\Gamma\nabla z,
\end{equation}
where $z$ is the height of the layer of sediment, $\alpha$ proportional to the coefficient  drag is taken to equal to $10^{-4}m^{-1}s^2,$ $1<\Gamma<3$ is a constant and $\gu$ stands for the velocity of the winds.
\item Bagnold flow is based on an energy approach and is modified by Gadd et al by including a critical speed for moving \cite{Bagnold1, Gadd}.
The flow $q_{s_0}$ is given by \begin{equation}
q_{s_0}=\frac{B_G}{\varphi_s}(\gu(z)-\gu_c(z))^3.
\end{equation}
assume $z = 1$ below the bottom, the hypothesus  logarithmic wire is satisfied.
\begin{equation}
\gu(z=1)=\frac{\sqrt{\tau_b}}{k\sqrt\tau}\ln\frac{7.5}{D_G}
\end{equation} which gives
\begin{equation}
q_{s_0}=\alpha(|\tau|_b^{\frac12}-\tau_{bc}^{\frac12})^3
\end{equation}
The total flow is then given by

\begin{equation}\label{eq3}q_s=q_{s_0}(\frac{\mathbf u}{\|\mathbf u\|}-\lambda_s\nabla z),\end{equation}

where $\mathbf u$ is the fluid velocity and $\lambda_s$ is the inverse value of the maximum slope
of the sediment surface when the velocity is 0.
\item The flow of  Komarova et al. \cite{komorova,komorova1} is given by
\begin{equation}\label{eq4}
q=\alpha|\tau|^{\frac12}\Big(\tau-\Gamma\nabla z\Big),\,\,
\end{equation}
with $\Gamma$ can take values between $\frac12$ and 6 and $\tau$ is given by $ \tau= \rho\frac{|\gu|^{2}}{C^{2}}\frac{\gu}{|\gu|}$ where $\gu$ is the velocity of the wind, $\rho$ the area density.
Komarova et al. \cite{komorova} consider that the speed is zero in bed bottom. The speed of  fluid in the bed bottom is not the parameter that moves the grains size but the wall friction.

\end{enumerate}
From these different flows of transport, we obtain different models describing the transport of the sand. Thus,
injecting (\ref{eq2}),\,(\ref{eq3}),\,(\ref{eq4}) in (\ref{eq001}), we get different formulas modeling the sand transport in the desert.
The first one  of sand transportation due to Gekerma is given by:
\begin{equation}\label{eq6}
\frac{\partial z}{\partial t}-\frac{\Gamma}{1-p}\Delta z=-\frac{\alpha}{1-p}\nabla\cdot(|\mathbf u|^3\frac{\mathbf u}{|\mathbf u|}),
\end{equation}
This partial differential equation is called the Gekerma's model.
The model due to Bagnold is given by
\begin{equation}\label{eq7}
\frac{\partial z}{\partial t}-\frac{\Gamma}{1-p}\nabla\cdot \Big(\alpha(|\tau_b|^{\frac12}-\tau_{bc}^{\frac12})^3( \frac{\mathbf u}{|\mathbf u|}-\lambda_s\nabla z)\Big)=0,
\end{equation}
and the  Komarova's one is given by
\begin{equation}\label{eq8}
\frac{\partial z}{\partial t}-\frac{\Gamma}{1-p}\nabla\cdot(\Gamma \alpha|\tau|^b\nabla z)=-\frac{1}{1-p}\nabla\cdot(\alpha|\tau|^b\tau).
\end{equation}
\subsection{Scaling and parameterized models}
In this section, we shall do the scaling techniques on the equations (\ref{eq6}, (\ref{eq7}),\,(\ref{eq8}) to write their dimensionless versions. Let us introduce a characteristic time $\bar t$ and a characteristic length $\bar L,$ then
\begin{equation}\label{eq9} t=t'\bar t,\,\,x=x'\bar L,\end{equation}
where $t'$ and $x'$  are the dimensionless variables.
Let us introduce also  the characteristic height of the dunes and velocity of wind by
\begin{equation}\label{eq10}
z(\bar t t', \bar L x')=z'(t,x)\bar z,\,\,\text{and}\,\,\mathbf u(\bar t t', \bar L x')=\mathbf u'(t,x)\bar u,
\end{equation}
where $z'$ and $\mathbf u'$ are the rescaled variables.\\
We shall, in the sequel give dimensionless models  already stated previously.\\
Using (\ref{eq9}) and (\ref{eq10}) and the fact that $\frac{\partial \bar z z'}{\partial t}=\frac{\bar z}{\bar t}\frac{\partial z'}{\partial t'}$ and $\Delta z=\frac{\bar z}{\bar L^2}\Delta z'$ we get for the Gekerma model's
\begin{equation}\label{kon1}
\frac{\partial z'}{\partial t'}-\frac{\bar t\Lambda}{(1-p)\bar L^2}\Delta z'=-\frac{\bar t}{\bar z}\frac{\alpha\bar  \gu^3}{(1-p)\bar L}\nabla'\cdot\Big(|\mathbf u'|^3\frac{\mathbf  u'}{|\mathbf u'|}\Big).
\end{equation}

The dimensionless Bagnold model's is given by

\begin{equation}\label{ban1}
\frac{\partial z'}{\partial t'}-\frac{\lambda\alpha\rho^{\frac32}\bar u^3\bar t}{\bar L^2(1-p)C^3}\nabla'\cdot\Big((\|\mathbf u'\|-\mathbf u_c)^3\nabla'z\Big)=
\frac{\alpha\rho^{\frac32}\bar u^3\bar t}{\bar L(1-p)C^3 \bar z}\cdot\nabla'\cdot
\Big((\|\mathbf u'\|-\mathbf u_c)^3\frac{\mathbf u'}{\|\mathbf u'\|}\Big)
\end{equation}

where  $C=\frac{1}{\kappa}\ln(\frac{30 \bar z}{D_G})$ with $D_G$ is the diameter of the sand grains.\\
The last model, due to Komarova is given by
\begin{equation}\label{lan1}
\frac{\partial z'}{\partial t'}-\frac{\bar t\Lambda|\bar\tau|^b}{1-p}\nabla'\cdot(|\tau '|^b\nabla z')=-\frac{\alpha\bar \tau |\bar\tau|^b\bar t}{\bar z(1-p)}\nabla'\cdot(|\tau'|^b\tau').
\end{equation}
Having these dimensionless models, we shall fix the characteristic values corresponding to different situations of dunes. We will consider  short, mean and long-term dynamics of the sand dunes of the desert. In this work, we consider various situations\\
\textbf{A})\textbf{A dynamical study of deformation of dunes in the short term}\\
Consider a small dune height $z =1m$ course  $ 300 m/year $. The mean residence time of the grains in these cases the dunes is worth $\bar t=100 days\sim 2400 hours\sim 8.6 10^6s.\,\,\bar t$  will be compared to the average period of day and night alternately winds which is of the order of $\frac{1}{\bar{w}}=13hours=4.7 10^4s$. We then get a small parameter $\epsilon$ given by $\epsilon=\frac{1}{\bar t\bar w}=\frac{1}{200}.$\\
In addition, the magnitude coefficient $\frac{\lambda}{1-p}$ is 1, then we get $\frac{\lambda}{1-p} = 1$ and $\frac{1}{1-p} = 2.$ \\
In the sequel, we are going to estimate the coefficients of the dimensionless models  (\ref{kon1}), (\ref{ban1}), (\ref{lan1}) to obtain their parameterized model.\\
\textbf{A}-1) Gerkerma's model is given by: \\
\begin{equation}\label{}
 \frac{\partial z^{'}}{\partial t^{'}} - \dfrac{\bar{t}\Lambda}{(1-p)\bar{L}^{2}}\Delta z^{'} = -\frac{\bar{t}}{\bar{z}} \dfrac{\alpha\bar{u}^{3}}{(1-p)\bar{L}}\nabla^{'}\cdot(|\gu^{'}|^3 \ \displaystyle\frac{\gu^{'}}{|\gu^{'}|}),
\end{equation}
 with $ \alpha=10^{-4} m^{-1} s^{2}$ and we take $\Lambda=3$,$\bar{u}=1m/s$ ,$\bar{L}=300m$ ,$\bar{z}=1m $ and  $\bar{t}=1 month.$ \\
 We first calculate the coefficients of the equation: \\
 $ \frac{1}{\bar{w}\bar{t}}\simeq \frac{1}{200}\simeq\varepsilon,$\\
 $ \dfrac{\bar{t}\Lambda}{(1-p)\bar{L}^{2}}\simeq\frac{3}{\varepsilon},$\\
 $\dfrac{\bar{t}\alpha\bar{u}^{3}}{\bar{z}(1-p)\bar{L}}\simeq 6.$\\
  So finally we get the following model: \\
  \begin{equation}\label{}
 \ \displaystyle\frac{\partial z^{'}}{\partial t^{'}} -\frac{3}{\varepsilon} \Delta z^{'} =
 6\nabla^{'}\cdot(|\gu^{'}|^3 \ \displaystyle\frac{\gu^{'}}{|\gu^{'}|}).
\end{equation}
\textbf{A}-2)The Komarova's model is given by:\\
 \begin{equation}\label{}
 \frac{\partial z}{\partial t} + \dfrac{1}{1-p}\nabla'\cdot(\alpha |\tau|^b \tau)  -\dfrac{1}{1-p}\nabla'\cdot( \Lambda\alpha |\tau|^b \nabla' f) = 0\ \ \ b = 1/2
\end{equation}
 with  $\Lambda =6$ and  $\alpha= 100.$ The shear stress $\tau$ is given by the relation $ \tau= \rho\frac{|\gu|^{2}}{C^{2}}\frac{\gu}{|\gu|}$
 with  $ \rho $ standing for the area density, $C=\frac{1}{k}ln(\frac{30\bar{z}}{D_{G}}$) with $ D_{G} $ being the sand grain diameter . After nondimensionalization of the equation we get: \\
 \begin{equation}\label{44z}
 \frac{\partial z^{'}}{\partial t^{'}} + \dfrac{\alpha\bar{u}^{3}\rho^{3/2}\bar{t}}{\bar{z}(1-p)C^{3}\bar{L}}\nabla'\cdot(|\gu'|^{2}\gu')  -\dfrac{\bar{u}\rho^{1/2}\bar{t}\Lambda\alpha}{(1-p)C^{3} \bar{L}^{2}}\nabla'\cdot(|\gu'|\nabla z^{'}) = 0
\end{equation}
  with $\Lambda=6$ and  $\alpha= 100.$
  In this case we take the sand grain diameter $ D_{G}=3\,10^{-4}$ and $k=0,4.$\\
So Can calculate the parameters in the equation:\\
$ \frac{1}{\bar{w}\bar{t}}\simeq \frac{1}{200}\simeq\frac{1}{\varepsilon},$\\
  $ C=\frac{1}{k}ln(\frac{30\bar{z}}{D_{G}})\simeq33,5,$\\
$\dfrac{\alpha\bar{u}^{3}\rho^{3/2}\bar{t}}{\bar{z}(1-p)C^{3}\bar{L}}\simeq\frac{3}{2\varepsilon},$\\
$ \dfrac{\bar{u}\rho^{1/2}\bar{t}\Lambda\alpha}{(1-p)C \bar{L}^{2}} \simeq\frac{3}{\varepsilon},$\\
with ( $\bar{u}=1m/s$ ,$\bar{L}=300m$ ,$\bar{z}=1m $ , $\bar{t}=100 days.$ \\
So finally we get the following model:\\
\begin{equation}\label{44z}
 \frac{\partial z'}{\partial t'} +\frac{3}{2\varepsilon}\nabla'\cdot(|\gu'|^{2}\gu')-\frac{3}{\varepsilon}\nabla'\cdot(|\gu'|\nabla z') = 0
\end{equation}
\textbf{B})\textbf{ Dynamic study of deformation of mean  term of dunes}\\
In the case of residence in medium term dunes,we consider a large dune height $z =10m $ travels  $30 m/years$.
   The mean residence time of the grains in these cases the dunes is 8 years so we have $\frac{1}{\bar{w}}=4 days$  the average period between two weather depressions, $\bar{L}=30 m$,$\bar{z}=10m$ , $\bar{t}=8years$ , $\bar{u}=1m/s$ $\frac{\lambda}{1-p} = 1$ and  $\frac{1}{1-p} = 2$.\\
\textbf{B}-1)The Gerkerma's model is given by: \\
\begin{equation}\label{}
 \ \displaystyle\frac{\partial z^{'}}{\partial t^{'}} - \dfrac{\bar{t}\Lambda}{(1-p)\bar{L}^{2}}\Delta z^{'} = -\frac{\bar{t}}{\bar{z}} \dfrac{\alpha\bar{u}^{3}}{(1-p)\bar{L}}\nabla^{'}\cdot(|\gu^{'}|^3 \ \displaystyle\frac{\gu^{'}}{|\gu^{'}|})
\end{equation}
 with $ \alpha=10^{-4} m^{-1} s^{2}$ and we take $\Lambda=3$.\\
 First we calculate the factors of the equation \\
 $ \frac{1}{\bar{w}\bar{t}}=\dfrac{4days}{8\times365days}\simeq 10^{-3}, $\\
 $ \dfrac{\bar{t}\Lambda}{(1-p)\bar{L}^{2}}\simeq\frac{16}{\varepsilon},$\\
  $\dfrac{\bar{t}\alpha\bar{u}^{3}}{\bar{z}(1-p)\bar{L}}\simeq\frac{1}{4\varepsilon}.$\\
  Then we have the following model:\\
  \begin{equation}\label{}
 \ \displaystyle\frac{\partial z'}{\partial t'} -\frac{16}{\varepsilon} \Delta z' = -\frac{1}{4\varepsilon}\nabla'\cdot(|\gu'|^3 \ \displaystyle\frac{\gu'}{|\gu'|}).
\end{equation}
\textbf{B}-2)The Komarova's model is given by:\\
\begin{equation}\label{44z}
 \frac{\partial z^{'}}{\partial t^{'}} + \dfrac{\alpha\bar{u}^{3}\rho^{3/2}\bar{t}}{\bar{z}(1-p)C^{3}\bar{L}}\nabla'\cdot(|\gu'|^{2}\gu')  -\dfrac{\bar{u}\rho^{1/2}\bar{t}\Lambda\alpha}{(1-p)C^{3} \bar{L}^{2}}\nabla'\cdot(|\gu'|\nabla' z^{'}) = 0
\end{equation}
 with $\Lambda =6$ ,$\alpha= 100$, $\bar{L}=100m$,$\bar{z}=10m$ , $\bar{t}=8years$ and $\bar{u}=1m/s$.\\
 In this case we take the sand grain diameter $ D_{G}=3\,10^{-4}$ and $k=0,4 $ . So we can calculate the parameters of the equation: \\
 $ \frac{1}{\bar{w}\bar{t}}=\dfrac{4days}{8\times365days}\simeq 10^{-3},$\\
  $ C=\frac{1}{k}ln(\frac{30\bar{z}}{D_{G}}=34,5,$\\
$\dfrac{\alpha\bar{u}^{3}\rho^{3/2}\bar{t}}{\bar{z}(1-p)C^{3}\bar{L}}=\frac{1}{10\varepsilon},$\\
$ \dfrac{\bar{u}\rho^{1/2}\bar{t}\Lambda\alpha}{(1-p)C\bar{L}^{2}}=\frac{39}{\varepsilon}, $\\
then, we have: \\
\begin{equation}\label{}
\frac{\partial z'}{\partial t'} +\frac{1}{10\varepsilon}\nabla'\cdot(|\gu'|^{2}\gu')  -\frac{39}{\varepsilon}\nabla'\cdot(|\gu'|\nabla' z') = 0
\end{equation}
\textbf{C})\textbf{Dynamic study of long-term deformation of dunes}\\
In the case of residence in medium term dunes,we consider a large dune height $\bar{z}=50m$ travels  $300m/years$. In these cases, the dunes is 2 centuries which will be compared to the average annual period is cyclical winds of about $\frac{1}{\bar{w}}=1 year$. Furthermore, the order of magnitude of coefficient $\frac{\lambda}{1-p}$ is 1 then we get $\frac{\lambda}{1-p} = 1$ and  $\frac{1}{1-p} = 2$ \\
\textbf{C}-1)The Gerkerma's model is given by: \\
\begin{equation}\label{}
 \ \displaystyle\frac{\partial z'}{\partial t'} - \dfrac{\bar{t}\Lambda}{(1-p)\bar{L}^{2}}\nabla' z^{'} = -\frac{\bar{t}}{\bar{z}} \dfrac{\alpha\bar{u}^{3}}{(1-p)\bar{L}}\nabla^{'}\cdot(|\gu^{'}|^3 \ \displaystyle\frac{\gu^{'}}{|\gu^{'}|}),
\end{equation}
 with $ \alpha=10^{-4} m^{-1} s^{2}$ and we take $\Lambda=3,\,\,\bar{u}=1m/s$ ,$\bar{L}=300m$  and  $\bar{z}=50m.$  \\
 We first calculate the factors of the equation \\
  $ \frac{1}{\bar{w}\bar{t}}\simeq \frac{1}{200}\simeq\varepsilon,$\\
 $ \dfrac{\bar{t}\Lambda}{(1-p)\bar{L}^{2}}\simeq\frac{5}{\varepsilon^{2}},$\\
  $\dfrac{\bar{t}\alpha\bar{u}^{3}}{\bar{z}(1-p)\bar{L}}\simeq\frac{1}{2\varepsilon}.$\\
  Then we have the following model:\\
  \begin{equation}\label{}
 \ \displaystyle\frac{\partial z'}{\partial t'} -\frac{5}{\varepsilon^{2}} \Delta z' =\frac{1}{2\varepsilon} \nabla'\cdot(|\gu'|^3 \ \displaystyle\frac{\gu'}{|\gu'|})
\end{equation}
\textbf{C}-2)The Komorava's model is given by:\\
\begin{equation}\label{44z}
 \frac{\partial z^{'}}{\partial t^{'}} + \dfrac{\alpha\bar{u}^{3}\rho^{3/2}\bar{t}}{\bar{z}(1-p)C^{3}\bar{L}}\nabla'\cdot(|\gu'|^{2}\gu')  -\dfrac{\bar{u}\rho^{1/2}\bar{t}\Lambda\alpha}{(1-p)C^{3} \bar{L}^{2}}\nabla'\cdot(|\gu'|\nabla' z^{'}) = 0
\end{equation}
  with $\Lambda =6,\,\,\alpha= 100,\,\,\bar{L}=300m,\,\,\bar{z}=50m,\,\,\bar{t}=200ans$ and  $\bar{u}=1m/s$.\\
We take the sand grain diameter $ D_{G}=3\,10^{-4}$  and  $k=0,4 $ .
  Then we have: \\
$ \frac{1}{\bar{w}\bar{t}}=\dfrac{1}{200ans}\simeq0,005\simeq\varepsilon,$\\
  $ C=\frac{1}{k}ln(\frac{30\bar{z}}{D_{G}}=39,$\\
$\dfrac{\alpha\bar{u}^{3}\rho^{3/2}\bar{t}}{\bar{z}(1-p)C^{3}\bar{L}}=\frac{9}{\varepsilon},$\\
$ \dfrac{\bar{u}\rho^{1/2}\bar{t}\Lambda\alpha}{(1-p)C\bar{L}^{2}}=\frac{5}{8\varepsilon^{2}}. $\\
By replacing the parameters in the equation we get:\\
\begin{equation}\label{44z}
 \frac{\partial z'}{\partial t'} +\frac{9}{\varepsilon}\nabla\cdot(|\gu'|^{2}\gu')  -\frac{5}{8\varepsilon^{2}}\nabla'\cdot(|\gu'|\nabla z') = 0.
\end{equation}

In all the following, we will remove the ' and consider that $z'(t,x)= z^\epsilon(t,x).$
\section{Existence results and uniqueness}
In this section, we are going to quote the different results of the paper. 
The models in the short, mean and long term dynamics of sandbanks in the desert can be rewritten as follows:
\begin{equation}\label{eq1}
\left\{\begin{array}{ccc}\frac{\partial z^\epsilon}{\partial t}-\frac{a}{\epsilon^j}\nabla\cdot\big(g^\epsilon(t,x)\nabla z^\epsilon\Big)=\frac{b}{\epsilon^i}\nabla\cdot f^\epsilon(t,x)
\,\,[0,T)\times \mathbb T^2,\\
z^\epsilon(0,x)=z_0(x),\,\,x\in\,\mathbb T^2,\end{array}\right.\end{equation}
where $i=0,1,\,j=0,1,2,\,\,z_0\in L^2(\torus^2),\,\,z_0\in L^2(\torus^2)$ and

\begin{equation}\label{0ff}
f^\epsilon(t,x)=\Big(g_c(|\mathbf u^\epsilon(t,x)|)\frac{\mathbf u^\epsilon(t,x)}{|\mathbf u^\epsilon(t,x)|}\Big),
\,\,\text{and}\,\,g^\epsilon(t,x)=g_a(|\mathbf u^\epsilon(t,x)|),\end{equation}

 where $g_a$ and $g_c$ satisfy the following hypotheses.

\begin{equation}\label{hyp1}\left\{ \begin{array}{c}
g_{a}\geq g_{c}\geq0,\,\,g_{c}(0)=g'_{c}(0)=0,\\
\exists\, d\geq0,~\sup_{u\in\mathbb{R}^{+}}|g_{a}(u)|+\sup_{u\in\mathbb{R}^{+}}|g'_{a}(u)|\leq d,~
    \sup_{u\in\mathbb{R}^{+}}|g_{c}(u)|+\sup_{u\in\mathbb{R}^{+}}|g'_{c}(u)|\leq d,\\
   \exists\, U_{thr}\geq0,\,\,\exists\, G_{thr}>0,\,\,\textrm{such that}\,\, u\geq U_{thr}\Longrightarrow g_{a}(u)\geq G_{thr}.\end{array}\right.
\end{equation}
Equation (\ref{eq1}) obtained in the models of sand transportation is similar to the one obtained in \cite{FaFreSe}, but (\ref{eq1}) is more general because it contains more cases.\\
The vectors fields $\gu^\epsilon :[0,T)\times\torus^2\rightarrow\mathbb R^2$ is the dimensionless wind velocity.
 \\
At first, we shall establish existence and uniqueness of $z^\epsilon(t,x)$ solution to (\ref{eq1}). And in the second case, we prove estimates in norms of $z^\epsilon,$ which do not depends on $\epsilon.$ 
\subsection{Existence and uniqueness for a weak solution}
Let us first suppose that

\begin{equation}\label{hyp2}
|f^\epsilon|,\,\,\big|\nabla\cdot f^\epsilon(t,x)\big|\leq \epsilon^ i\gamma,\,i=0,1,
\end{equation}
where $\gamma$ is a constant not depending on $\epsilon.$
\begin{theorem}\label{theo1}
Let $\epsilon>0,\,a>0$ and $b$ and $c$ reals. Under hypotheses (\ref{hyp1}) and (\ref{hyp2}) and if $z_0\in H^1(\mathbb T^2),$ for any $T>0,$ there exists a unique solution $z^\epsilon\in L^2((0,T],L^2(\mathbb T^2))$ to (\ref{eq1}) satisfying
\begin{equation}\label{ada1.1}\int_{\mathbb T^2}\frac{\partial z^\epsilon}{\partial t}dx=0,\end{equation} and
\begin{equation}
\|z^\epsilon\|_{L^2((0,T)\times\mathbb T^2)}\leq \gamma.
\end{equation}
\end{theorem}
Before proving the above main existence result we need to show preliminaries key results.
As the diffusion coefficient of equation (\ref{eq1}) may cancel, we will start by regularizing equation (\ref{eq1}). Let $\nu>0,$ we consider the following regularizing equation:
 \begin{equation}\label{eq1nu}
\left\{\begin{array}{ccc}\frac{\partial z^{\epsilon,\nu}}{\partial t}-\frac{a}{\epsilon^j}\nabla\cdot\big((g^\epsilon(t,x)+\nu)\nabla z^{\epsilon,\nu}\Big)=\frac{b}{\epsilon^i}\big(f^\epsilon(t,x)\big)
\,\,\text{in}\,\,[0,T)\times \torus^2,\\
z^{\epsilon,\nu}(x,0)=z_0(x)\,\,\text{in}\,\,\torus^2.\end{array}\right.\end{equation}
They are stated and showed as follows.
\begin{prop1} Under the same hypotheses as in theorem\,\,\ref{theo1}, for all $\epsilon>0,$ for all $\nu>0$ their exists a unique solution $z^{\epsilon,\nu}\in L^2([0,T), L^2(\torus^2))$ to (\ref{eq1nu}).

\end{prop1}
\begin{proof}
Let $z^{\epsilon,\nu}_1$ and $z^{\epsilon,\nu}_2$ two solutions of (\ref{eq1nu}). Their difference is  solution to
\begin{equation}\label{eq1nu1}
\left\{\begin{array}{ccc}\frac{\partial (z_1^{\epsilon,\nu}-z_2^{\epsilon,\nu})}{\partial t}-\frac{a}{\epsilon^j}\nabla\cdot\big((g^\epsilon(t,x)+\nu)\nabla( z_1^{\epsilon,\nu}-z_2^{\epsilon,\nu})\Big)=0\,\,\text{in}\,\,[0,T)\times \torus^2,\\
z_1^{\epsilon,\nu}(0,x)-z_2^{\epsilon,\nu}(0,x)=0\,\,\text{in}\,\,\torus^2,\end{array}\right.\end{equation}
Multiplying (\ref{eq1nu1}) by $(z_1^{\epsilon,\nu}-z_2^{\epsilon,\nu})$ and integrating over $\mathbb T^2$ we get

$$\frac12\frac{d}{dt}\|z_1^{\epsilon,\nu}-z_2^{\epsilon,\nu}\|_2^2+\frac{a}{\epsilon^j}\int_{\torus^2}(g^\epsilon+\nu)\nabla |z_1^{\epsilon,\nu}-z_2^{\epsilon,\nu}|^2dx=0$$
which gives because of the positivity of the second term of the equality
\begin{equation}
\frac12\frac{d}{dt}\|z_1^{\epsilon,\nu}-z_2^{\epsilon,\nu}\|_2^2\leq0.
\end{equation}
Integrating this last equality from 0 to $t$, we get
\begin{equation}
\|z_1^{\epsilon,\nu}(t)-z_2^{\epsilon,\nu}(t)\|_2^2\leq0.
\end{equation}
From this last equality, we get $z_1^{\epsilon,\nu}(t)=z_2^{\epsilon,\nu}(t),\,\,\forall t\in [0,T),$ which gives uniqueness.\\
Existence of $z^{\epsilon,\nu}$ follows from adaptations of results of Ladyzenskaja, Solonnikov and Ural'Ceva \cite{Lady} and Lions \cite{Lions}.
In the following, we prove estimates of $z^{\epsilon,\nu}$ solution to (\ref{eq1nu}) which does not depend on $\nu$ and $\epsilon.$
\begin{prop1}\label{prop1} Under the same hypotheses as in the theorem \ref{theo1}, $\forall \nu>0$ the sequence of solutions $z^{\epsilon,\nu}$ solution to (\ref{eq1nu}) satisfies the following inequalities:
\begin{equation}\label{esti1}
\|z^{\epsilon,\nu}\|_{L^\infty((0,T),L^2(\mathbb T^2))}\leq \sqrt{\|z_0\|_2^2+b\tilde\gamma T},
\end{equation}
\begin{equation}\label{esti2}
\|\nabla z^{\epsilon,\nu}\|^2_{L^\infty((0,T),L^2(\mathbb T^2))}\leq \frac{\epsilon^j}{aG_{thr}}\Big(\|z_0\|_2^2+2b\tilde\gamma\Big),\,\,j=0,1.
\end{equation}
\begin{equation}\label{esti3}
\Big\|\frac{\partial z^{\epsilon,\nu}}{\partial t}\Big\|_{L^2([0,T), L^2(\mathbb T^2))}\leq \epsilon^j\frac{b\gamma}{G_{thr}},\,\,j=0,1.
\end{equation}

\end{prop1}
\begin{proof}
For all $\epsilon>0,\,\,i=0,1,\,\,j=0,1,2$ multiplying $(\ref{eq1nu})$ by $z^{\epsilon,\nu}$ and integrating over $\mathbb T^2,$ we obtain
\begin{equation}\label{grand1}\frac12\frac{d\|z^{\epsilon,\nu}\|_2^2}{dt}+\int_{\mathbb T^2}\frac{a}{\epsilon^j}((g^\epsilon+\nu)|\nabla z^{\epsilon,\nu}|^2dx\leq \int_{\mathbb T^2}\frac{b}{\epsilon^i}| f^\epsilon||\nabla z^{\epsilon,\nu}|dx.\end{equation}
Integrating (\ref{grand1}) from 0 to $T$ we get
\begin{equation}
\frac12\|z^{\epsilon,\nu}(T)\|_2^2+\int_0^T\int_{\mathbb T^2}\frac{a}{\epsilon^j}((g^\epsilon+\nu)|\nabla z^{\epsilon,\nu}|^2dxdt\leq \frac{b\gamma}{\epsilon^i}\int_0^T\Big(\int_{\mathbb T^2}|\nabla z^{\epsilon,\nu}|^2dx\Big)^{\frac12}dt+\frac12\|z_0\|^2
\end{equation}
From this last inequality, using (\ref{hyp1}), we get:
\begin{equation}G_{thr}\int_0^T\int_{\mathbb T^2}|\nabla z^{\epsilon,\nu}|^2\leq \frac{b\gamma\epsilon^{j-i}}{a}\int_0^T\Big(\int_{\mathbb T^2}|\nabla z^{\epsilon,\nu}|dx\Big)^{\frac12}dt+\frac{\epsilon^j}{2a}\|z_0\|^2.
\end{equation}
There exists a constant  $\tilde\gamma$ depending only  on $\|z_0\|^2,\,\,G_{thr},\,\,a$ and $b$  such that
\begin{equation} \label{eer}\|\nabla z^{\epsilon,\nu}\|_{L^2([0,T), L^2(\torus^2))}\leq \tilde\gamma;
\end{equation}
Using the fact that $(g^\epsilon+\nu)>0$ and $a>0$ we obtain from (\ref{grand1})
$$\frac12\frac{d\|z^{\epsilon,\nu}\|_2^2}{dt}\leq \int_{\mathbb T^2}\frac{b}{\epsilon^i}| f^\epsilon||\nabla \cdot z^{\epsilon,\nu}|dx$$
and using hypotheses (\ref{hyp2}) and (\ref{eer}), we get
\begin{equation}
\frac{d\|z^{\epsilon,\nu}\|_2^2}{dt}\leq b\tilde\gamma.
\end{equation}
Integrating this last inequality, from 0 to $t\in[0, T)$ we get
\begin{equation}
\|z^{\epsilon,\nu}(t)\|_2\leq \sqrt{\|z_0\|_2^2+b\tilde\gamma T},
\end{equation} and then  $\sup_{t\in[0,T)}\|z^{\epsilon,\nu}(t)\|_2\leq\sqrt{\|z_0\|_2^2+b\tilde\gamma T}.$
This last inequality gives (\ref{esti1}).\\
From (\ref{grand1}) we get
\begin{eqnarray}
{\frac12\int_{\{t\in [0,T):\,\,|\mathbf u|\geq U_{thr}\}}\frac{d\|z^{\epsilon,\nu}\|_2^2}{dt}dt+\int_{\{t\in [0,T):\,\,|\mathbf u| \geq U_{thr}\}}\int_{\mathbb T^2}\frac{a}{\epsilon^j}((g^\epsilon+\nu)|\nabla z^{\epsilon,\nu}|^2dxdt {} }~~~~~~
\nonumber\\
\leq \int_{0}^T\int_{\mathbb T^2}\frac{b}{\epsilon^i}|\nabla \cdot f^\epsilon||z^{\epsilon,\nu}|dxdt.\end{eqnarray}
From this last inequality, using (\ref{hyp1}), we have
\begin{eqnarray}{G_{thr}\frac{a}{\epsilon^j}\int_{\{t\in [0,T):\,\,|\mathbf u| \geq U_{thr}\}}\int_{\mathbb T^2}|\nabla z^{\epsilon,\nu}|^2dxdt {} }~~~~~~
\nonumber\\
\leq \int_{0}^T\int_{\torus^2}\frac{b}{\epsilon^i}|\nabla \cdot(f^\epsilon)||z^{\epsilon,\nu}|dxdt+2\|z_0\|_2^2+b\tilde\gamma T.\end{eqnarray}
There exists $t_0\in \{t\in[0,T):\,\,|\mathbf u|<U_{thr}\}$ such that
\begin{equation}
G_{thr}\int_{\mathbb T^2}|\nabla z^{\epsilon,\nu}(t_0)|^2dx\leq\frac{\epsilon^j}{a}\Big(2b\tilde\gamma+2\|z_0\|_2^2\Big),\end{equation} which gives

\begin{equation}\|\nabla z^{\epsilon,\nu}(t_0,\cdot)\|_2^2\leq \frac{\epsilon^j}{aG_{thr}}\Big( 2b\tilde\gamma+2\|z_0\|_2^2\Big),\,\,j=0,1.\end{equation}
Multiplying (\ref{eq1nu}) by $\displaystyle\frac{\partial z^{\epsilon,\nu}}{\partial t}$ and integrating over $\mathbb T^2,$ we get
\begin{equation}\int_{\mathbb T^2}\Big|\frac{\partial z^{\epsilon,\nu}}{\partial t}\Big|^2dx-\frac{a}{\epsilon^j}\int_{\mathbb T^2}\nabla\cdot\Big((g^\epsilon+\nu)\nabla z^{\epsilon,\nu}\Big)\frac{\partial z^{\epsilon,\nu}}{\partial t}dx=\frac{b}{\epsilon^i}\int_{\mathbb T^2}(\nabla\cdot f^\epsilon)\frac{\partial z^{\epsilon,\nu}}{\partial t}dx.\end{equation}

\begin{equation}
\int_{\mathbb T^2}\Big|\frac{\partial z^{\epsilon,\nu}}{\partial t}\Big|^2dx+\frac{a}{\epsilon^j}\int_{\mathbb T^2}(g^\epsilon+\nu)\nabla z^{\epsilon,\nu}\nabla\big(\frac{d z^{\epsilon,\nu}}{dt}\big)dx=\frac{b}{\epsilon^i}\int_{\mathbb T^2}(\nabla\cdot f^\epsilon)\frac{\partial z^{\epsilon,\nu}}{\partial t}dx,
\end{equation}
which gives
\begin{equation}
\int_{\mathbb T^2}\Big|\frac{\partial z^{\epsilon,\nu}}{\partial t}\Big|^2dx+\frac{a}{2\epsilon^j}\int_{\mathbb T^2}(g^\epsilon+\nu)\Big|\nabla \frac{\partial z^{\epsilon,\nu}}{\partial t}\Big|^2dx\leq\frac{b}{\epsilon^i}\int_{\mathbb T^2}|f^\epsilon|\Big|\nabla\frac{\partial z^{\epsilon,\nu}}{\partial t}\Big|dx,
\end{equation}
using the fact that the first term is positive, we get

\begin{equation}
\frac{a}{\epsilon^j}\Big\|\sqrt{g^\epsilon+\nu}|\nabla\frac{\partial z^{\epsilon,\nu}}{\partial t}|\Big\|_2^2\leq b\gamma\Big\|\nabla\frac{\partial z^{\epsilon,\nu}}{\partial t}\Big\|_2
\end{equation}
which integrating in $\{t\in[0,T),\,\,|\mathbf u|< U_{thr}\},$ gives
\begin{equation}
\frac{a}{\epsilon^j}G_{thr}\int_{t\in[0,T),\,\,|\mathbf u|< Uthr}\Big\|\nabla\frac{\partial z^{\epsilon,\nu}}{\partial t}\Big\|^2_2dt\leq b\gamma\int_{t\in[0,T),\,\,|\mathbf u|< Uthr}\Big\|\nabla\frac{\partial z^{\epsilon,\nu}}{\partial t}\Big\|_2dt.
\end{equation}
From this last inequality, there exists $t_0\in \{t\in[0,T),\,\,|\mathbf u|< Uthr\}$ such that
\begin{equation}\label{zt}
\Big\|\nabla\frac{\partial z^{\epsilon,\nu}}{\partial t}(t_0,\cdot)\Big\|_2\leq\epsilon^j\frac{b\gamma}{aG_{thr}}\,\,\text{for a given}\,\,t_0.
\end{equation}
From (\ref{zt}), we get using Fourier series of $\frac{\partial z^{\epsilon,\nu}}{\partial t}$ that

\begin{equation}
\Big\|\frac{\partial z^{\epsilon,\nu}}{\partial t}(t_0)\Big\|_2\leq \epsilon^j\frac{b\gamma}{aG_{thr}}.
\end{equation}
Then \begin{equation}
\sup_{t_0\in[0,T)}\Big\|\frac{\partial z^{\epsilon,\nu}}{\partial t}(t_0)\Big\|_2\leq \epsilon^j\frac{b\gamma}{aG_{thr}}.
\end{equation}

which gives, estimate (\ref{esti3}).
\end{proof}

\end{proof}

\begin{proof}\,\, of \textbf{Theorem \ref{theo1}}
Existence of $z^\epsilon$ solution to $(\ref{eq1})$ follows from proposition \ref{prop1}. As estimates (\ref{esti1}),(\ref{esti2}),(\ref{esti3}) do not depends on $\nu,$ letting $\nu\rightarrow0,$ we obtain $z^\epsilon$ with the same properties.\\
Integrating (\ref{eq1}) with respect to $x$  over $\torus^2$, we get equality (\ref{ada1.1}).
\end{proof}
\subsection{Existence in a general case}
In the following, we shall consider the general case i.e. we dot not assume (\ref{hyp2}) for $i=1$. From (\ref{hyp1}) we prove that  $\nabla\cdot f^\epsilon$ is in fact bounded by a constant $\gamma.$
Since we aim to get qualitative information about the  asymptotic behavior of $z^\epsilon$ as $\epsilon$ goes to 0, we need estimates of $z^\epsilon$ which do not depends on $\epsilon$ or which is bounded when $\epsilon$ goes to zeros.
To prove estimates of $z^\epsilon$ which do not depends on $\epsilon$ we are going to consider  the following periodic case.

As the frequency of winds is considered periodic in desert, we can assume that
\begin{equation}\label{ff0}\left\{\begin{array}{ccc}\mathbf u^\epsilon(t,x)=\mathcal U(t,\frac t\epsilon,x)\\
\textbf{and}\,\, \theta\longmapsto \mathcal U(t,\theta,x)\,\,\textbf{is periodic of period 1}.\end{array}\right.\end{equation}

 Then
$\theta\longmapsto g_a(|\mathbf u(t,\theta,x)|),\,\,g_c(|\mathbf u(t,\theta,x)|)\frac{\mathbf u(t,\theta,x)}{|\mathbf u(t,\theta,x)|}$ where $g_a$ and $g_c$ satisfy the hypotheses (\ref{hyp1}) are also periodic of period 1, and let us set
\begin{equation}\label{ff1}
g^\epsilon(t,x)=g_\epsilon(t,\frac t\epsilon,x)\,\,\text{and}\,\, f^\epsilon(t,x)=f_\epsilon(t,\frac t\epsilon,x),\end{equation}
where $g^\epsilon$ and  $f^\epsilon$ are difined in (\ref{0ff}).

\begin{theorem}\label{theo2}
Let $\epsilon>0,\,a>0$ and $b$ and $c$ reals. Under hypotheses  (\ref{0ff})and (\ref{hyp1}) and if $z_0\in H^1(\mathbb T^2),$ for any $T>0,$ there exists a unique solution $z^\epsilon\in L^2((0,T]\times\mathbb T^2)$ to (\ref{eq1}) satisfying
\begin{equation}\label{ada1}\int_{\mathbb T^2}\frac{\partial z^\epsilon}{\partial t}dx=0,\end{equation} and
\begin{equation}
\|z^\epsilon\|_{L^2((0,T)\times\mathbb T^2)}\leq \gamma
\end{equation}
for a constant $\gamma$ not depending on $\epsilon.$
\end{theorem}

The proof theorem \ref{theo2} is inspired by Faye et al  works in \cite{FaFreSe} where they study the long-term dynamics of sand dunes in areas subjected to the tide.

In the following we are going to focus our efforts on existence and uniqueness of time-space periodic parabolic equations.
From this, we then get existence of the solution to equation (\ref{eq1}). Existence of $z^\epsilon$ in a time interval
depending on $\epsilon$ is a straightforward adaptation of results from LadyzensKaja,
Solonnikov and Ural' Ceva \cite{Lady} or Lions \cite{Lions}. Our aim is to prove that $z^\epsilon$ solution to (\ref{eq1}) is bounded
independently of $\epsilon.$ And let us introduce the following regularized equations stated as follows:

 find $Z^\nu=Z^\nu(t,\theta,x)$ periodic of period 1 in $\theta$ solution to
\begin{equation}\label{eqreg12}
\frac{\partial Z^\nu}{\partial \theta}-\frac{a}{\epsilon^j}\nabla\cdot\Big((g_\epsilon+\nu)\nabla Z^\nu\Big)=\nabla\cdot f_\epsilon\,\,\text{in}
\,\,[0,T)\times \mathbb T^2,\,\,j=0,1.\end{equation}
From hypotheses (\ref{hyp1}) and (\ref{ff1}), functions $f_\epsilon$ and $g_\epsilon$ satisfy the following hypotheses

\begin{equation}\label{hyp02}\left\{\begin{array}{ccc}
\theta\mapsto ( g_\epsilon,  f_\epsilon)\,\,\text{is periodic of period 1},\\
x\mapsto ( g_\epsilon, f_\epsilon)\,\,\text{is defined on}\,\, \mathbb T^2\\
| g_\epsilon|\leq \gamma,\,\,| f_\epsilon|\leq \gamma,\,\, \Big|\frac{\partial  f_\epsilon}{\partial t}\Big|\leq \epsilon^2d,\,\,\Big|\frac{\partial g_\epsilon}{\partial t}\Big|\leq \epsilon d,\\
\Big|\frac{\partial f_\epsilon}{\partial \theta}\Big|\leq \epsilon d,\,\,\Big|\frac{\partial g_\epsilon}{\partial \theta}\Big|\leq \epsilon d\\
\Big|\nabla g_\epsilon\Big|\leq \epsilon d,\,\,\Big|\nabla \cdot f_\epsilon\Big|\leq \epsilon d, \Big|\frac{\partial \nabla\cdot f_\epsilon}{\partial t}\Big|\leq \epsilon^2 d,\,\,\,| g_\epsilon|\leq d f_\epsilon.
\end{array}\right.\end{equation}

\begin{equation}\label{ff2}\left\{\begin{array}{ccc}
\exists \tilde G_{thr}\geq0,\,\,\theta_\alpha,\,\theta_\omega\in[0,1]
\\
 \tilde G_{thr}\leq  g_\epsilon
\end{array}\right.\end{equation}
\begin{theorem}\label{thnu}
Under hypotheses (\ref{hyp1}),(\ref{ff0}),\ref{ff1}, (\ref{hyp02}), (\ref{ff2}), $\forall \epsilon>0$ and $\nu>0,$ there exists a unique solution $Z^\nu\in L^\infty_{\#}(\mathbb R, L^2(\torus^2))$ solution to (\ref{eqreg12}). \\
Moreover this solution satisfies the following inequalities
\begin{equation}\label{th1}
\int_{\theta_\alpha}^{\theta_\omega}\int_{\mathbb T^2}|\nabla Z^\nu|^2dx\,d\theta\leq \epsilon^i\frac{\gamma}{G_{thr}},\,\,i=0,\,1,
\end{equation}
\begin{equation}\label{th2}
\|\nabla Z^\nu(\theta_0,\cdot)\|_2\leq \frac{\gamma\epsilon^{i}}{\sqrt{\tilde G_{thr}}},\,\,i=0,1
\end{equation}
\begin{equation}\label{th3}
\| Z^\nu(\theta_0,\cdot)\|_{L_{\#}^\infty(\mathbb R, L^2(\torus^2))}\leq \frac{\gamma\epsilon^{i}}{\sqrt{\tilde G_{thr}}}+2\gamma.
\end{equation}
\begin{equation}\label{th4}
\| \frac{\partial Z^\nu}{\partial t}(\theta_0,\cdot)\|_{L_{\#}^\infty(\mathbb R, L^2(\torus^2))}\leq \frac{\gamma f(\epsilon)}{\sqrt{\tilde G_{thr}}},
\end{equation}
where $f(\epsilon)\rightarrow \lambda>0$ as $\epsilon\rightarrow0.$
\end{theorem}
\begin{proof}
Multiplying (\ref{eqreg12}) by $S^\nu$ and integrating in $\mathbb T^2$ and in $[0,1],$
 we get
 \begin{equation}
 \int_0^1\int_{\mathbb T^2}\frac{\partial Z^\nu}{\partial \theta} Z^\nu dx\,d\theta+\frac{a}{\epsilon ^i}\int_0^1\int_{\mathbb T^2}(g_\epsilon+\nu)|\nabla Z^\nu|^2dxd\theta\leq \int_0^1\int_{\mathbb T^2}|f_\epsilon||\nabla Z^\nu|dx\,d\theta
 \end{equation}
 From this inequality, we get
 \begin{equation}\frac{a}{\epsilon ^i}\int_0^1\int_{\mathbb T^2}g_\epsilon|\nabla Z^\nu|^2dxd\theta\leq \int_0^1\int_{\mathbb T^2}|f_\epsilon||\nabla Z^\nu|dx\,d\theta
 \end{equation}
 then from (\ref{ff2}) we have \begin{equation}\label{h3}
 \|\sqrt{\tilde G_{thr}}|\nabla Z^\nu|\|_{L^2(\mathbb R, L^2(\mathbb T^2))}\leq \frac{\epsilon^i\gamma}{a}.
 \end{equation}
 Integrating from $\theta_\alpha$ to  $\theta_\omega$ we get
 \begin{equation}\label{eqreg1} \int_{\theta_\alpha}^{\theta_\omega}\int_{\mathbb T^2}|\nabla Z^\nu|^2dx\,d\theta\leq \epsilon^i\frac{\gamma}{a\sqrt{\tilde G_{thr}}},\,\,i=0,\,1.
 \end{equation}
 Using (\ref{eqreg1}), there exists $\theta_0\in [\theta_\alpha,\theta_\omega]$ such that
 \begin{equation}\|\nabla Z^\nu(\theta_0,\cdot)\|_2^2\leq \epsilon^i\frac{\gamma}{a\sqrt{\tilde G_{thr}} },\,\,i=0,\,1. \end{equation}
 Using Fourier series of $Z^\nu,$ we can prove that
 \begin{equation}\label{th4}\|Z^\nu(\theta_0,\cdot)\|_2^2\leq\|\nabla Z^\nu(\theta_0,\cdot)\|_2^2\leq \epsilon^i\frac{\gamma}{a\sqrt{\tilde G_{thr}}},\,\,i=0,\,1. \end{equation}
  Multiplying (\ref{eqreg12}) by $Z^\nu,$ and integrating over $\mathbb T^2,$ we get

  \begin{equation}\label{eqreg23}\frac12\frac{d\|Z^\epsilon(\theta,\cdot)\|_2^2}{d\theta}+\frac{a}{\epsilon^i}\int_{\{x\in\mathbb T^2,\,\,g_a(|\mathcal U(t,\theta,x)|)=0\}}\nu|\nabla Z^\nu|^2dx+\frac{a}{\epsilon^i}\int_{\{x\in\mathbb T^2,\,\,g_a(|\mathcal U(t,\theta,x)|)\neq0\}}(g_\epsilon+\nu)|\nabla Z^\nu|^2dx$$
  $$\leq \int_{\{x\in\mathbb T^2,\,\,g_a(|\mathcal U(t,\theta,x)|)\neq0\}}|f^\epsilon\cdot\nabla Z^\nu|dx\end{equation}
  But we have,
  \begin{equation}\label{eqreg24}
  \int_{\{x\in\torus^2,\,\,g_a(|\mathcal U(t,\theta,x)|)\neq0\}}\big|f^\epsilon\cdot\nabla Z^\nu\big|dx\leq\int_{\{x\in\torus^2,\,\,g_a(|\mathcal U(t,\theta,x)|)\neq0\}}\frac{(g_\epsilon+\nu)}{4}|\nabla Z^\nu|^2dx$$
  $$+\int_{\{x\in\torus^2,\,\,g_a(|\mathcal U(t,\theta,x)|)\neq0\}} \frac{|f^\epsilon|^2}{g^\epsilon+\nu}dx.
  \end{equation}

  Using (3.50) in (3.49) and passing the first term in the right hand side in the left hand side, we get
  \begin{equation}\frac{d\|Z^\nu(\theta,\cdot)\|^2}{d\theta}\leq 2\gamma.
  \end{equation}
  Integrating this last inequality from $\theta_0$ to  another $\theta\in[\theta_0,\theta_0+1]$ we get using (\ref{th4}) the following
  \begin{equation}\|Z^\nu(\theta,\cdot)\|_2^2\leq \epsilon^i\frac{\gamma}{a\sqrt{\tilde G_{thr}}}+2\gamma,\,\,i=0,\,1.
  \end{equation}
  Since $Z^\nu$ is periodic, we get inequality (\ref{th3}).
  \\
  $\ds\frac{\partial Z^\nu}{\partial t}$ is solution to
  \begin{equation} \label{tt}
  \frac{\partial \Big(\frac{\partial Z^\nu}{\partial t}\Big)}{\partial \theta}-\frac{a}{\epsilon^{j}}\nabla\cdot\Big((g_\epsilon+\nu)\nabla\Big( \frac{\partial Z^\nu}{\partial t}\Big)\Big)=\frac{a}{\epsilon^j}\nabla \Big(\frac{\partial g_\epsilon}{\partial t}\nabla Z^\nu\Big)+\nabla\cdot\Big(\frac{\partial f_\epsilon}{\partial t}\Big),
  \end{equation} and  multiplying (\ref{tt}) by $\ds\frac{\partial Z^\nu}{\partial t}$ and integrating in $x\in \torus^2,$ we get
  \begin{equation}\label{r5}
  \frac12\frac{\partial\Big\|\frac{\partial Z^\nu}{\partial t}\Big\|_2^2}{\partial \theta}+\frac{a}{\epsilon^{j}}\int_{\torus^2}(g_\epsilon+\nu)\Big|\nabla \frac{\partial Z^\nu}{\partial t}\Big|^2dx\leq \int_{\torus^2}\Big|\frac{\partial f_\epsilon}{\partial t}\Big|\Big|\nabla\frac{\partial Z^\nu}{\partial t}\Big|dx+\frac{a}{\epsilon^{j}}\int_{\torus^2}\Big|\frac{\partial g_\epsilon}{\partial t}\Big||\nabla Z^\nu|\Big|\nabla\frac{\partial Z^\nu}{\partial t}\Big|dx.
  \end{equation}
   But the first term of right hand side can be written using hypothesis (\ref{hyp02})and (\ref{hyp1}) as follows
 \begin{equation}\label{hhh6}
   \int_{\torus^2}\Big|\frac{\partial f_\epsilon}{\partial t}\Big|\Big|\nabla\frac{\partial Z^\nu}{\partial t}\Big|dx\leq \gamma \int_{\torus^2}\sqrt{g_\epsilon}\Big|\nabla\frac{\partial Z^\nu}{\partial t}\Big|dx$$
  $$\leq \gamma\Big\|\sqrt{g_\epsilon}\Big|\nabla\frac{\partial Z^\nu}{\partial t}\Big\|_2.
  \end{equation}
  From hypothesis (\ref{hyp02}) we deduce
\begin{equation}\label{hhh7}
  \int_{\torus^2}\Big|\frac{\partial g_\epsilon}{\partial t}\Big|
 |\nabla Z^\nu|\Big|\nabla\frac{\partial Z^\nu}{\partial t}\Big|dx
 \leq
 \Big\|\sqrt{\Big|\frac{\partial g_\epsilon}{\partial t}\Big|}|\nabla Z^\nu|\Big\|_2\Big\|\sqrt{\Big|\frac{\partial g_\epsilon}{\partial t}\Big|}
 \Big|\nabla\frac{\partial Z^\nu}{\partial t}
 \Big|\Big\|_2$$
 $$\leq \gamma^2\Big\|\sqrt{ g_\epsilon}|\nabla Z^\nu|\Big\|_2\Big\|\sqrt{g_\epsilon}\Big|\nabla\frac{\partial Z^\nu}{\partial t}\Big|\Big\|_2.
  \end{equation}
 From this inequality , using (3.55) and (3.56) we have \\
   \begin{equation}
  \frac12\frac{\partial\Big\|\frac{\partial Z^\nu}{\partial t}\Big\|_2^2}{\partial \theta}+\frac{a}{\epsilon^{j}}\int_{\torus^2}(g_\epsilon+\nu)\Big|\nabla \frac{\partial Z^\nu}{\partial t}\Big|^2dx\leq  \gamma\Big\|\sqrt{g_\epsilon}\Big|\nabla\frac{\partial Z^\nu}{\partial t}\Big\|_2+  \frac{a}{\epsilon^{j}}\gamma^2\Big\|\sqrt{ g_\epsilon}|\nabla Z^\nu|\Big\|_2\Big\|\sqrt{g_\epsilon}\Big|\nabla\frac{\partial Z^\nu}{\partial t}\Big|\Big\|_2.
  \end{equation}
  Integrating from 0 to  1 we get\\
   \begin{equation}
  \frac{a}{\epsilon^{j}}\int_{0}^{1}\int_{\torus^2} g_\epsilon \Big|\nabla \frac{\partial Z^\nu}{\partial t}\Big|^2dxd\theta\leq \int_{0}^{1} \gamma\Big\|\sqrt{g_\epsilon}\Big|\nabla\frac{\partial Z^\nu}{\partial t}\Big\|_2d\theta+  \int_{0}^{1}\frac{a}{\epsilon^{j}}\gamma^2\Big\|\sqrt{ g_\epsilon}|\nabla Z^\nu|\Big\|_2\Big\|\sqrt{g_\epsilon}\Big|\nabla\frac{\partial Z^\nu}{\partial t}\Big|\Big\|_2d\theta.
  \end{equation}
  Using Holder inequality in the second term,we get \\
    \begin{equation}
  \frac{a}{\epsilon^{j}}\Big\|\sqrt{g_\epsilon }\Big|\nabla \frac{\partial Z^{\nu}}{\partial t}\Big\|_{L^2(\mathbb R, L^2(\mathbb T^2))}^{2}\leq  \gamma\Big\|\sqrt{g_\epsilon}\Big|\nabla\frac{\partial Z^\nu}{\partial t}\Big|\Big\|_2+  \frac{a}{\epsilon^{j}}\gamma^2\Big\|\sqrt{ g_\epsilon}|\nabla Z^\nu|\Big\|_{L^2(\mathbb R, L^2(\mathbb T^2))}\Big\|\sqrt{g_\epsilon}\Big|\nabla\frac{\partial Z^\nu}{\partial t}\Big|\Big\|_{L^2(\mathbb R, L^2(\mathbb T^2))}.
  \end{equation}
  which gives\\
   \begin{equation}
  \frac{a}{\epsilon^{j}}\Big\|\sqrt{g_\epsilon } \Big|\nabla \frac{\partial Z^{\nu}}{\partial t}\Big|\Big\|_{L^2(\mathbb R, L^2(\mathbb T^2))}\leq  \gamma +  \frac{a}{\epsilon^{j}}\gamma^2\Big\|\sqrt{ g_\epsilon}|\nabla Z^\nu|\Big\|_{L^2(\mathbb R, L^2(\mathbb T^2))}.
  \end{equation}
 Using (\ref{h3}) we have\\
     \begin{equation}
  \Big\|\sqrt{g_\epsilon } \Big|\nabla \frac{\partial Z^{\nu}}{\partial t}\Big|\Big\|_{L^2(\mathbb R, L^2(\mathbb T^2))}\leq  \frac{\epsilon^{j}}{a}\gamma +  \epsilon^{j}\gamma^3.
  \end{equation}
   From hypothesis (\ref{hyp02}) we get\\
   \begin{equation}
  \Big\|\Big|\nabla \frac{\partial Z^{\nu}}{\partial t}\Big\|_{L^2(\mathbb R, L^2(\mathbb T^2))}\leq\gamma\frac{f(\epsilon)}{\sqrt{\tilde G_{thr}}} .
  \end{equation}
  Integrating from $\theta_\alpha$ to  $\theta_\omega$ we get\\
  \begin{equation}\label{r1}
 \int_{\theta_\alpha}^{\theta_\omega} \Big\|\nabla \frac{\partial Z^{\nu}}{\partial t}\Big\|_{L^2(\mathbb R, L^2(\mathbb T^2))}d\theta\leq\gamma\frac{f(\epsilon)}{\sqrt{\tilde G_{thr}}} .
  \end{equation}
  From (\ref{r1}), there exists $\theta_0\in [\theta_\alpha,\theta_\omega]$ such that\\
  \begin{equation}\label{r2}
 \Big\|\nabla \frac{\partial Z^{\nu}(\theta_0,.)}{\partial t}\Big\|_{L^2(\mathbb R, L^2(\mathbb T^2))}\leq\gamma\frac{f(\epsilon)}{\sqrt{\tilde G_{thr}}} .
  \end{equation}
  Using Fourier series of $Z^\nu,$ we can prove that\\
   \begin{equation}\label{r3}
 \Big\|\Big|\frac{\partial Z^{\nu}(\theta_0,.)}{\partial t}\Big|\Big\|_{L^2(\mathbb R, L^2(\mathbb T^2))}\leq\Big\|\Big|\nabla \frac{\partial Z^{\nu}(\theta_0,.)}{\partial t}\Big|\Big\|_{L^2(\mathbb R, L^2(\mathbb T^2))}\leq\gamma\frac{f(\epsilon)}{\sqrt{\tilde G_{thr}}} ,
  \end{equation}
  which gives (\ref{th4}).
  \end{proof}
  \begin{theorem}
 Under hypothesis (\ref{hyp1}),(\ref{ff0}),(\ref{ff1}), (\ref{hyp02}) and(\ref{ff2})
 there exists $Z=Z(t,\theta,x)\in L^\infty_{\#}(\mathbb R, L^2(\mathbb T^2))$ , solution to \\
 \begin{equation}\label{eqreg122}
\frac{\partial Z}{\partial \theta}-\frac{a}{\epsilon^j}\nabla\cdot\Big(g_\epsilon\nabla Z\Big)=\nabla\cdot f_\epsilon\,\,\text{in}\,\,
\,\,[0,T)\times \mathbb T^2.\end{equation}
Moreover,this solution satisfies the following inequalities:
\begin{equation}\label{th10}
\int_{\theta_\alpha}^{\theta_\omega}\int_{\mathbb T^2}|\nabla Z|^2dx\,d\theta\leq \epsilon^i\frac{\gamma}{G_{thr}},\,\,i=0,\,1,
\end{equation}
\begin{equation}\label{th20}
\|\nabla Z(\theta_0,\cdot)\|_2\leq \frac{\gamma\epsilon^{i}}{\sqrt{\tilde G_{thr}}},\,\,i=0,1
\end{equation}
\begin{equation}\label{th30}
\| Z(\theta_0,\cdot)\|_{L_{\#}^\infty(\mathbb R, L^2(\torus^2))}\leq \frac{\gamma\epsilon^{i}}{\sqrt{\tilde G_{thr}}}+2\gamma.
\end{equation}
\begin{equation}\label{th40}
\| \frac{\partial Z}{\partial t}(\theta_0,\cdot)\|_{L_{\#}^\infty(\mathbb R, L^2(\torus^2))}\leq \frac{\gamma f(\epsilon)}{\sqrt{\tilde G_{thr}}},
\end{equation}
where $f(\epsilon)\rightarrow \lambda>0$ as $\epsilon\rightarrow0.$

  \end{theorem}
 \begin{proof} Existence of $Z$ follows from theorem \ref{thnu}. In fact it suffices that to let $\nu$ going to zero and the desired result is obtained letting $\nu$ to zero.
  \end{proof} 
%
%
\begin{proof} of \textbf{Theorem \ref{theo2}}
  \\ Having this theorem on hand, let $Z^\epsilon(t,x)= Z(t,\frac t\epsilon,x).$ Then $z^\epsilon(t,x)- Z(t,\frac t\epsilon,x)$ is solution to

  \begin{equation}\left\{\begin{array}{ccc}\label{eqZ}
  \frac{\partial(z^\epsilon- Z^\epsilon)}{\partial t}-\frac{a}{\epsilon^j}\nabla\cdot\Big(g^\epsilon\nabla (z^\epsilon- Z^\epsilon)\Big)= \frac{\partial Z}{\partial t}(t,\frac t\epsilon,x)
  \\
  z^\epsilon(0,x)-Z^\epsilon(0,x)=z_0(x)-Z(0,0,x);\end{array}\right.\end{equation}
  Multiplying (\ref{eqZ}) by $z^\epsilon- Z^\epsilon$ and integrating over $\torus^2$, we get
  \begin{equation}
  \|z^\epsilon(t,\cdot)-Z^\epsilon(t,\cdot)\|_2\leq \gamma.
  \end{equation}
  From this last inequalities, we get that $z^\epsilon(t,x)$ is not far from $Z^\epsilon(t,x)$ and as $Z^\epsilon(t,x)=Z(t,\frac{t}{\epsilon},x)$
  is bounded in $L_{\#}^\infty(\mathbb R, L^2(\torus^2))$ we conclude that $z^\epsilon$ is also bounded in $L_{\#}^\infty(\mathbb R, L^2(\torus^2)).$

\end{proof}
\section{Homogenization results}
Let us consider (\ref{eq1})with coefficients given by (\ref{0ff}). We are interested by the behavior of $z^\epsilon$ when $\epsilon\rightarrow0.$ In fact our aim in this section is to study the homogenization problem associated to (\ref{eq1}). By theorem \ref{theo1}, we have showed that, there exists a unique solution $z^\epsilon$ to (\ref{eq1}) which is bounded independently in $\epsilon$ in $L^\infty(\mathbb R, L^2(\mathbb T^2)).$\\
It is obvious that,
\begin{equation}\label{hom1}
g^\epsilon(t,x)\,\,\text{two scale converges to}\,\,\tilde g(t,\theta,x)\in L^\infty([0,T), L^\infty(\mathbb R, L^2(\torus^2)))
\end{equation}
and

\begin{equation}\label{hom2}
f^\epsilon(t,x)\,\,\text{two scale converges to}\,\,\tilde f(t,\theta,x)\in L^\infty([0,T), L^\infty(\mathbb R, L^2(\torus^2)))
\end{equation}

where $\tilde g$ and $\tilde f$ are given by

\begin{equation}\label{gg}
\tilde g(t,\theta,x)=g_a(|\mathcal U(t,\theta,x)|)\,\,\text{and}\,\,\tilde g(t,\theta,x)=g_a(|\mathcal U(t,\theta,x)|)\frac{\mathcal U(t,\theta,x)}{|\mathcal U(t,\theta,x)|}.
\end{equation}
\begin{theorem}
Under assumptions (\ref{ff0}),\,\, (\ref{ff1}) and (\ref{hyp02})-(\ref{ff2}), for any $T^,$ not depending on $\epsilon,$ the sequence $(z^\epsilon)$  of solutions to (\ref{eq1}), with coefficients given by (\ref{0ff}), two-scale converges to the profile $U\in L^\infty([0,T), L^\infty_{\#}(\mathbb R, L^2(\mathbb R, L^(\torus^2)))$ which is characterized by:
\begin{enumerate}
  \item \begin{equation}\label{eq1180}
  \frac{\partial U}{\partial\theta} - \nabla\cdot\big( \tilde{g}\nabla U) = \nabla\cdot \tilde f
  \end{equation}
  in the short and mean terms of dynamics of dunes
 and
 \item
  \begin{equation}\label{eq11200}
  \nabla\cdot\big( \tilde{g}\nabla U) = 0
  \end{equation} in the case of long terms dynamics of dunes.
\end{enumerate}
where  $\tilde g$ et $\tilde f$ are given by (\ref{gg}).
\end{theorem}
\begin{proof}
Defining test function $\psi^{\epsilon}(t,x)=\psi(t,\frac{t}{\epsilon},x)$ for any $\psi(t,\theta,x)$,regular with
compact support in $ \,\,[0,T)\times \mathbb T^2,$ and periodic in $ \theta $ with period 1. Then multiplying (\ref{eq1}) by $\psi^{\epsilon}$ and integrating in $ \,\,[0,T)\times \mathbb T^2,$ we have \\
\begin{equation}\label{eq111}
\int_{\mathbb T^2}\int^{T}_{0}\frac{\partial z^\epsilon}{\partial t}\psi^{\epsilon}dtdx -\frac{a}{\epsilon^j} \int_{\mathbb T^2}\int^{T}_{0}\nabla\cdot\big(g^\epsilon(t,x)\nabla z^\epsilon\Big)\psi^{\epsilon}dtdx =\frac{b}{\epsilon^i}\int_{\mathbb T^2}\int^{T}_{0}\nabla\cdot f^\epsilon(t,x)\psi^{\epsilon}dtdx \end{equation}
Then integrating by parts in the first integral in $[0,T)$ and using the Green formula in $\mathbb T^2$ in the second integral we have
\begin{equation}\label{eq112}
-\int_{\mathbb T^2} z_{0}^\epsilon(x)\psi(0,0,x)dx -\int_{\mathbb T^2}\int^{T}_{0}\frac{\partial \psi^\epsilon}{\partial t}z^{\epsilon}dtdx  +\frac{a}{\epsilon^j} \int_{\mathbb T^2}\int^{T}_{0} g^\epsilon(t,x)\nabla z^\epsilon \nabla\psi^{\epsilon}dtdx $$
$$=\frac{b}{\epsilon^i}\int_{\mathbb T^2}\int^{T}_{0}\nabla\cdot f^\epsilon(t,x)\psi^{\epsilon}dtdx \end{equation}
Again thanks to the Green formula in the third integral we get \\
\begin{equation}\label{eq113}
-\int_{\mathbb T^2} z_{0}^\epsilon(x)\psi(0,0,x)dx-\int_{\mathbb T^2}\int^{T}_{0}\frac{\partial \psi^\epsilon}{\partial t}z^{\epsilon}dtdx  -\frac{a}{\epsilon^j} \int_{\mathbb T^2}\int^{T}_{0} z^\epsilon \nabla\cdot\big( g^\epsilon(t,x)\nabla\psi^{\epsilon})dtdx $$
$$=\frac{b}{\epsilon^i}\int_{\mathbb T^2}\int^{T}_{0}\nabla\cdot f^\epsilon(t,x)\psi^{\epsilon}dtdx \end{equation}
But $$\frac{\partial\psi^\epsilon}{\partial t}=\big(\frac{\partial\psi}{\partial t}\big)^{\epsilon}+ \dfrac{1}{\epsilon} \big(\frac{\partial\psi}{\partial \theta}\big)^{\epsilon}.$$
And then we have\\
 \begin{equation}\label{eq114}
  -\int_{\mathbb T^2}\int^{T}_{0} \Big((\frac{\partial\psi}{\partial t})^{\epsilon}+ \dfrac{1}{\epsilon} \big(\frac{\partial\psi}{\partial \theta}\big)^{\epsilon} + \frac{a}{\epsilon^j}\nabla\cdot\big( g^\epsilon(t,x)\nabla\psi^{\epsilon})\Big)z^\epsilon dtdx$$
   $$=\frac{b}{\epsilon^i}\int_{\mathbb T^2}\int^{T}_{0}\nabla\cdot f^\epsilon(t,x)\psi^{\epsilon}dtdx + \int_{\mathbb T^2} z_{0}^\epsilon(x)\psi(0,0,x)dx \end{equation}
  Using the two-scale convergence due to Nguetseng \cite{Ngue} and Allaire \cite{allaire:1992}, if a sequence $ f^{\epsilon}$ is bounded in $L^\infty(0,T, L^2(\torus^2))$,then there exists a profile $ U(t,\theta,x)$, periodic of period 1 with respect to
  $\theta $ , such that for all $ \psi(t,\theta,x)$, regular with compact support with respect to $ (t,x)$ and periodic of period 1 with respect to $\theta$ ,we have \\
  \begin{equation}\label{eq1106}
   \int_{\mathbb T^2}\int^{T}_{0} f^\epsilon \psi^{\epsilon}dtdx\longrightarrow \int_{\mathbb T^2}\int^{T}_{0} \int^{1}_{0} U \psi dtdxd\theta
   \end{equation}
 In the case of short and mean term dynamics of dunes, the exponent of the parameter $\epsilon$ in the equation  is $i=j=1$  or $j=1$ and $i=0.$\\
Multiplying (4.9) by $\epsilon$ 
and using (\ref{eq1106}) we have\\
   \begin{equation}\label{eq115}
  \int_{\torus^2}\int^{T}_{0} \Big[\epsilon\Big(\frac{\partial\psi}{\partial t}\Big)^{\epsilon}+\Big(\frac{\partial\psi}{\partial \theta}\Big)^{\epsilon} + \nabla\cdot\big( g^\epsilon(t,x)\nabla\psi^{\epsilon}\big)\Big]z^\epsilon dtdx =$$
   $$\int_{\torus^2}\int^{T}_{0} \nabla\cdot f^\epsilon\psi^\epsilon dxdt
  +\epsilon\int_{\torus^2} z_0^\epsilon(x)\psi(0,0,x)dx
  \end{equation}
  As $g^\epsilon$ and $f^\epsilon$ are bounded(see(\ref{hyp02})and $\psi^\epsilon$ is a regular function,
  $g^\epsilon(t,x)\nabla\psi^{\epsilon}$ and $\nabla\psi^{\epsilon}$ can be considered as test functions. Using (\ref{hom1}) and (\ref{hom2})  we have the equation satisfied by $U :$\\
  \begin{equation}\label{eq118}
  \frac{\partial U}{\partial\theta} + \nabla\cdot\big( \tilde{g}\nabla U) = \nabla \tilde f.
  \end{equation}
  In the case of long term dynamics of dunes, i.e. $j=i=2$ or $j=2$ and $i=1,$ then multiplying (4.9) by $\epsilon^{2}$ and using (\ref{eq1106}) we have\\
    \begin{equation}\label{eq119}
 \int_{\torus^2}\int^{T}_{0} \Big[\epsilon^2\Big(\frac{\partial\psi}{\partial t}\Big)^{\epsilon}+\epsilon\Big(\frac{\partial\psi}{\partial \theta}\Big)^{\epsilon} + \lim_{\epsilon\longrightarrow 0}\int_{\mathbb T^2}\int^{T}_{0} \big(\nabla\cdot\big( g^\epsilon(t,x)\nabla\psi^{\epsilon})\big)z^\epsilon dtdx $$
 $$=\int_{\mathbb T^2}\int^{T}_{0}\nabla\cdot f^\epsilon\psi^\epsilon dxdt+\epsilon^2\int_{\torus^2}z_0(x)dx
  \end{equation}
 Passing to the limit as $\epsilon\longrightarrow 0,$  using two scale convergence 
 we obtain from (4.13) the equation satisfied by the two-scale limits $U:$\\
  \begin{equation}\label{eq1120}
  \nabla\cdot\big( \tilde{g}\nabla U) = 0.
  \end{equation}
\end{proof}
\section{A corrector result}
The two-scales convergence result shows that the solution $z^\epsilon$ of  equation (\ref{eq1}) can be expressed as follows:
$$z^\epsilon(t,x)=\sum_{i=0}^{+\infty}\epsilon^iU^{i}(t,\frac t\epsilon,x),$$ where $\theta\longmapsto U(t,\theta,x)$ is periodic of period 1. The aim in this section is to characterize the equation satisfied by $U^1.$
Following the idea of two-scales convergence, we suppose that the coefficients $g^\epsilon$ and $f^\epsilon$ two scale converge respectively to $\tilde g$ and $\tilde f.$ This leads to writing

\begin{equation}\label{cor1}g^\epsilon(t,x)=\tilde g^\epsilon(t,x)\end{equation} and
\begin{equation}\label{cor2}f^\epsilon(t,x)=\tilde f^\epsilon(t,x)\end{equation} where
\begin{equation}\label{cor3}\tilde g^\epsilon(t,x)=\tilde g(t,\frac t\epsilon,x)\,\,\text{and}\,\,\tilde f^\epsilon(t,x)=\tilde f(t,\frac t\epsilon,x)\end{equation} where $\tilde g$ and $\tilde f$ are the two scales limits of $g^\epsilon$ and $f^\epsilon.$


\begin{theorem}
Under assumptions (\ref{ff0}),\,\, \ref{ff1},  and (\ref{hyp02})-(\ref{ff2}), for any $T^,$ not depending on $\epsilon,$ considering the sequence $(z^\epsilon)$  of solutions to (\ref{eq1}), with coefficients given by (\ref{cor1}),\,\,(\ref{cor2}) and $U^\epsilon(t,x)=U(t,\frac t\epsilon,x)$ where $U$ is solution to (\ref{eq1180}),
the following estimate holds for $z^\epsilon-U^\epsilon,$
\begin{equation}
\Big\|\frac{z^\epsilon-U^\epsilon}{\epsilon}\Big\|_{L^\infty([0,T), L^2(\mathbb T^2))}\leq \alpha.
\end{equation}

Furthermore, the sequence $(\frac{z^\epsilon-U^\epsilon}{\epsilon})$ two-scale converge to the profile \\$U^1\in L^\infty([0,T), L^\infty_{\#}(\mathbb R, L^2(\mathbb R, L^(\torus^2)))$ which is the unique solution to
 \begin{equation}\label{cor4}
 \frac{\partial U^1}{\partial \theta}-\nabla\cdot\Big(\tilde g\nabla U^1\Big)=\frac{\partial U}{\partial t}.
 \end{equation}

\end{theorem}
\begin{proof}
We have
\begin{equation} \frac{\partial U^\epsilon}{\partial t}=\Big( \frac{\partial U}{\partial t}\Big)^\epsilon+\frac{1}{\epsilon}\Big( \frac{\partial U}{\partial \theta}\Big)^\epsilon,
\end{equation}
where \begin{equation}
\Big( \frac{\partial U}{\partial t}\Big)^\epsilon(t,x)=\frac{\partial U}{\partial \theta}(t,\frac t\epsilon,x)\,\,\text{and}\,\, \Big( \frac{\partial U}{\partial \theta}\Big)^\epsilon=\frac{\partial U}{\partial \theta}(t,\frac t\epsilon, x)\end{equation}
then  $U^\epsilon$ is solution to
\begin{equation}
\frac{\partial U^\epsilon}{\partial t}-\frac1\epsilon\nabla\cdot\Big(\tilde g^\epsilon\nabla U^\epsilon\Big)=\frac{1}{\epsilon}\nabla\cdot\tilde f^\epsilon+\Big( \frac{\partial U}{\partial t}\Big)^\epsilon.\end{equation}
Using (\ref{cor1}) and (\ref{cor2}) and equation (\ref{eq1}), one can proof that  $z^\epsilon-U^\epsilon$ is solution to

\begin{equation}
\frac{\partial\Big(\frac{ z^\epsilon-U^\epsilon}{\epsilon}\Big)}{\partial t}-\frac1\epsilon\nabla\cdot\Big(\tilde g^\epsilon\nabla\big(\frac{z^\epsilon-U^\epsilon}{\epsilon}\big)\Big)=\frac1\epsilon\Big(\frac{\partial U}{\partial t}\Big)^\epsilon.
\end{equation}
Following the idea developed in \cite{FaFreSe} and a result of Ladyzenskaja, Sollonikov and  Ural'Ceva \cite{Lady}, $(\frac{\partial U}{\partial t}\Big)^\epsilon$ is solution to a parabolic linear and bounded coefficient, then is bounded. Then, using the same
arguments as in the proof of theorem 1.1 we obtain that $\frac{ z^\epsilon-U^\epsilon}{\epsilon}$ is bounded, and it two-scale
 converges to a profile
$U^1\in L^\infty([0, T];L^\infty_{\#}(\mathbb R; L^2(\mathbb T^2)))$ and that this
profile $U^1$ satisfies equation (\ref{cor4}).
\end{proof}


\begin{thebibliography}{99}

\bibitem{allaire:1992}
G.~Allaire, \emph{{H}omogenization and Two-Scale convergence}, SIAM J.
  Math. Anal. \textbf{23}, 1482--1518, 1992.
  \bibitem{AndClauDou}  B.~Andreotti, P. Claudin, S.~Douady, \emph{Selection of dunes shpes and velocities. Part 1: Dynamics of sand, wind and barchans}, Eur. Phys J. B., p341-352, 2002.
      \bibitem{Bagnold1} R.A. Bagnold, \emph{ The physics of blown sand and desert dunes}, Chapmann and Hall,1941.
\bibitem{Bagnold1} R. A. Bagnold, \emph{The movement of desert sand}, Proceedings of the Royal
  Society of London A \textbf{157}, 594--620, 936.
\bibitem{CokeWarreGou} R.Coke, A. Warren, \emph{Desert Geomorphology}, UCL Press, 1993.
\bibitem{Hersen} P. Hersen,\emph{Morphog\`enese et dynamique des Barchanes}, Th\`ese de doctorat, 2004.


\bibitem{FaFreSe} I.~ Faye, E.~ Fr\'enod, D.~ Seck, \emph{Singularly perturbed degenerated parabolic equations and application to seabed
 morphodynamics in tided environment},Discrete and Continuous Dynamical Systems, Vol 29 $N^{o} 3$ March 2011, pp 1001-1030.
\bibitem{Gadd}  P. E. Gadd, W. Lavelle and D. J. P. Swift,\emph{Estimates of sand transport on the New York shelf
using near-bottom current meter observations}, J. Sed. Petrol., 48, 239-252, 1978.



\bibitem{Gek} T.~ Gekerma,\emph{A linear stability analysis of tidally generated sands waves}, J. Fluid Mech, vol 417, 303-322, 2000.


 \bibitem{komorova} N. L. Komarova and A. Newell, \emph{Nonlinear dynamics of sandbanks and
sandwaves}. J. Fluid Mech., vol. 415, pp. 285-321, 2000.
\bibitem{komorova1} N. L. Komarova and S. J. M. Hulcher, \emph{Linear mechanisms for sand wave formation}. J. Fluid Mech.,  413, pp. 219-246, 2000.


 \bibitem{com1} Morsi et Alexander, A. J.  An investigation of Particle Trajectories inn Two Fluid Rech., 251, 661 - 685, 1972.
 \bibitem{Lady} O. A. Ladyzenskaja, V. A. Solonnikov and N. N. Ural'ceva, \emph{Linear and Quasi-Linear Equations
of Parabolic Type}, AMS Translation of Mathematical Monographs 23, 1968.

\bibitem{Lions} J. L. Lions, \emph{Remarques sur les \'equations diff\'erentielles ordinaires}, Osaka Math. J., 15, 131-142, 1963.

 \bibitem{com2} H. Nishimori,  and N. Ouchi, \emph{Formation of Ripple Patterns and Dunes by Wind}, Blown Sand. Phys. Rev. Lett., \textbf{71,} 197-200, 1993.
 \bibitem{Ngue} G. Nguetseng, A general convergence result for a functional related to the theory of homogenization, SIAM J. Math. Anal. 20 , 608–623, 1989.
\bibitem{com3}B. B. Willets, and M. A.  RICE \emph{Collisions of quartz grains with a sand bed: the influence of incident angle} , Earth Surface Processes and Landforms ,14,719-730, 1989.
\end{thebibliography}
\end{document}